\documentclass{amsart}
\usepackage{amsmath}
\usepackage{amsfonts}
\usepackage{amssymb}
\usepackage{amsthm}
\usepackage{newlfont}
\usepackage{graphicx}
\usepackage{longtable}
\usepackage{tabularx}
\usepackage{array, ragged2e}
\newtheorem{thm}{Theorem}
\newtheorem{cor}{Corollary}
\newtheorem{lem}{Lemma}
\newtheorem{prop}{Proposition}
\theoremstyle{definition}
\newtheorem{defn}{Definition}
\newenvironment{acknowledge}{\noindent\textbf{Acknowledgments.}}{}
\numberwithin{equation}{section} \numberwithin{thm}{section}
\numberwithin{cor}{section} \numberwithin{lem}{section}
\numberwithin{prop}{section} \numberwithin{defn}{section}
\newcommand{\set}[1]{\left\{{#1}\right\}}

\newcommand{\coord}[1]{\left({#1}\right)}

\newcommand{\intpart}[1]{\left\lfloor{#1}\right\rfloor}
\newcommand{\fracpart}[1]{\left\{{#1}\right\}}
\newcommand{\roof}[1]{\left\lceil{#1}\right\rceil}
\newcommand{\coprime}[1]{\left(#1\right)}
\newcommand{\N}{\mathbb{N}}
\newcommand{\Z}{\mathbb{Z}}
\newcommand{\Q}{\mathbb{Q}}

\newcommand{\sconv}[1]{\mathrm{conv}\set{#1}}

\newcommand{\lfrac}[1]{\frac{\lambda_{#1}\kappa}{h}}
\newcommand{\modsum}[1]{\lambda_{#1}\kappa\text{ (mod }h\text{)}}
\newcommand{\modpart}[1]{\fracpart{\frac{\lambda_{#1}\kappa}{h}}}
\newcommand{\invmodpart}[1]{\fracpart{\frac{\lambda_{#1}(h-\kappa)}{h}}}
\newcommand{\wholepart}[1]{\intpart{\frac{\lambda_{#1}\kappa}{h}}}
\newcommand{\rfrac}[1]{\roof{\frac{\kappa\lambda_{#1}}{h}}}

\newcommand{\con}[1]{c_{#1}^{(\kappa)}}
\newcommand{\ai}[2]{a_{#1}^{(#2)}}

\begin{document}
\title{Toric Fano $3$-folds with Terminal Singularities}
\author{Alexander M Kasprzyk}
\thanks{Mathematical
Sciences, University of Bath, Bath, BA2 7AY, United Kingdom.
\newline\indent e-mail~\texttt{A.M.Kasprzyk@maths.bath.ac.uk}.}
\begin{abstract}
This paper classifies all toric Fano $3$-folds with terminal
singularities. This is achieved by solving the equivalent
combinatoric problem; that of finding, up to the action of
$GL(3,\Z)$, all convex polytopes in $\Z^3$ which contain the
origin as the only non-vertex lattice point.
\end{abstract}
\maketitle
\section{Background and Introduction}\label{intro}
A \emph{toric variety} of dimension $n$ over an algebraically
closed field $k$ is a normal variety $X$ that contains a torus
$T\cong (k^*)^n$ as a dense open subset, together with an action
$T\times X\rightarrow X$ of $T$ on $X$ that extends the natural
action of $T$ on itself.

Let $M\cong\Z^n$ be the lattice of characters of the torus $T$,
with dual lattice $N=Hom(M,\Z)$. Then every toric variety $X$ has
an associated fan $\Delta$ in $N\otimes\Q$. The converse also
holds; to any fan $\Delta$ there is an associated toric variety
$X(\Delta)$. For details on the construction and deeper properties
of this crucial bijection see~\cite{Dan78,Oda78,Ful93,CCAG,Cox02}.

A normal variety $X$ is a \emph{Fano variety} if some multiple of
the anticanonical divisor $-K_X$ is an ample Cartier divisor. Thus
$X$ is $\Q$-Gorenstein. There is a well known description of what
it means for a toric variety $X$ to be a Fano variety in terms of
its fan $\Delta$: Let $\set{\rho_i}_{i\in I}$ be the set of rays
of $\Delta$. For each $i\in I$ there exists a unique primitive
lattice element of $\rho_i$, which by a traditional abuse of
notation we continue to denote $\rho_i$. Then $X$ is Fano if and
only if $\set{\rho_i}_{i\in I}$ correspond to the vertices of a
convex polytope in $N$ (see~\cite{Dan78,CCAG}).

Fano varieties are important for several reasons. In particular
they play a significant role in the Minimal Model Program
(see~\cite{Wis01,FS03}). Many general results concerning smooth
toric Fano $n$-folds are known~\cite{Wis01}; of particular
relevance, it has been shown that there are precisely $18$ smooth
toric Fano $3$-folds.

A normal variety $X$ is \emph{$\Q$-factorial} if every prime
divisor $\Gamma\subset X$ has a positive integer multiple
$c\Gamma$ which is a Cartier divisor. Once again, for the toric
case there exists a well known description in terms of the fan.
The toric variety $X$ is $\Q$-factorial if and only if the fan
$\Delta$ is simplicial (see~\cite{Oda78,Dais02}).

We say that a fan $\Delta$ is \emph{terminal} if each cone
$\sigma\in\Delta$ satisfies the following:

\begin{itemize}
 \item[(i)]the rays $\rho_1,\ldots,\rho_k$ of $\sigma$ are
 contained in an affine hyperplane $H:(m(n)=1)$ for some $m\in
 M_{\Q}$
 \item[(ii)]there are no other elements of the lattice $N$ in the
 part of $\sigma$ under or on $H$
 $$\text{ie. } N\cap\sigma\cap
 (m(n)\leq1)=\set{0,\rho_1,\ldots,\rho_k}.$$
\end{itemize}

A toric variety $X$ is \emph{terminal} (ie. has terminal
singularities~\cite{Reid83M}) if and only if the fan $\Delta$ is
terminal.

Terminal singularities play an important role in birational
geometry. A great deal of classification results exist for various
cases; for example the results
of~\cite{Mori82,Reid83M,Mor84,Mori88}. In~\cite{Mor85} it was
shown that, with two exceptions, isolated canonical cyclic
quotient singularities in dimension three are all either
Gorenstein or terminal. In~\cite{Reid85} Reid addresses the issue
of classifying $3$-fold terminal singularities. In the notation
of~\cite[Theorem 6.1,\ II]{Reid85} we are in the case
$\frac{1}{r}\coord{a,-a,1,0;0}$.

We are now in a position to state our aim, a complete
classification of all toric Fano $3$-folds with terminal
singularities, in terms of a purely combinatorial problem. Namely
we wish to find, up to the action of $GL(3,\Z)$, all convex
lattice polytopes in $\Z^3$ which contain only the origin as a
non-vertex lattice point (by which we mean that no lattice points
lie on the surface of the polytope other than the vertices, and no
lattice points are contained in the interior of the polytope other
than the origin).

An equivalent restatement for $2$-folds can be found in~\cite[pp.
192-3]{CCAG}; precisely five polytopes are found, of which two are
minimal (the Fano triangle and the Fano square, which make an
appearance in Section~\ref{classifyMinNow}) and one is maximal, in
the sense of Definitions~\ref{minimalDef} and~\ref{maximalDef}.
The approach used for this classification relies on the basic
result that, up to the action of $GL(2,\Z)$, there is a unique
lattice point free triangle (namely $\sconv{0,e_1,e_2}$). This
fails to hold for three dimensions (see~\cite{Sca85}). It is also
worth observing that in dimension two all polytopes are simplicial
(and hence the corresponding toric variety is, at worst,
$\Q$-factorial), something which is clearly not the case in
dimension three.

The classification presented in this paper is inspired by the work
of A.Borisov and L.Borisov~\cite{BB,BB92}. Results given
in~\cite{BB92,Bor00} assure us that a finite classification is
possible. The combinatoric approach we adopt is based on that
formulated in~\cite{BB}. In this unpublished work, the essential
steps described can be outlined thus:

\begin{itemize}
 \item[(i)]Observe that every polytope can be ``grown'' from a ``minimal'' polytope.
 \item[(ii)]These minimal polytopes divide into
 tetrahedra and non-tetrahedra.
 \item[(iii)]The minimal tetrahedra can be classified in terms of their
 barycentric coordinates.
 \item[(iv)]The minimal non-tetrahedra can be determined directly.
 \item[(v)]A recursive algorithm can be written, allowing a computer to ``grow'' these
 minimal polytopes and hence classify all polytopes of interest.
\end{itemize}

The result of Proposition~\ref{Tet2Sum} is a specific case
of~\cite[Proposition 3]{BB92}, however the proof presented here is
of an elementary combinatorial nature, in keeping with the style
of the remainder of this paper. In addition the results of
Table~\ref{tetTable} are obtained more explicitly than
in~\cite{BB92}; again the justification for repeating these
results lies in the methods used to obtain them. With a nice
restatement of Proposition~\ref{Tet2Sum} (concerning tetrahedra
containing one non-vertex lattice point) we obtain a result which
closely mirrors one of~\cite{Sca85} (concerning tetrahedra
containing no non-vertex lattice points), although once more the
methods of proof are very different.

For practical reasons the final classification is not contained in
this paper, but has been made available on the internet (see the
end of Section~\ref{finalClassAll} for the address). We conclude
this introduction by presenting a summary of the main features of
this list (see also Table~\ref{3DNumTable}).
\begin{thm}\label{mainSummary}
Up to isomorphism, there are exactly $233$ toric Fano $3$-folds
having at worst $\Q$-factorial singularities, of which $18$ are
smooth. There are an additional $401$ having terminal
singularities that are not $\Q$-factorial.

There are $12$ minimal cases with at worse $\Q$-factorial singularities: eight with Picard number $1$, two
with Picard number $2$, and two with Picard number $3$. There is one minimal case with terminal singularities, corresponding to a polytope with $5$ vertices.

There are nine maximal cases, corresponding to polytopes with $8$ (three
occurrences), $9$, $10$ (two occurrences), $11$, $12$ and $14$ vertices. Only
those with $8$ vertices are are $\Q$-factorial.
\end{thm}
\begin{acknowledge}
The author would like to express his gratitude to Professor
Alexander Borisov for making~\cite{BB} available; the current
paper relies heavily upon many of the ideas developed in this
unpublished work. Special thanks to Dr. Gregory Sankaran for
introducing me to this problem, and for his invaluable
explanations and advice.

The author wishes to acknowledge funding by a Mathematics CASE
studentship from EPSRC/GCHQ.
\end{acknowledge}
\section{Classifying the Tetrahedra: The Barycentric Coordinates}
We will make frequent appeals to the following well known result:
\begin{lem}\label{emptyTriEqu}
Any lattice point free triangle with vertices
$\set{0,x_1,x_2}\subset\Z^3$ is equivalent under the action of
$GL(3,\Z)$ to the triangle with vertices $\set{0,e_1,e_2}$.
\end{lem}
Let $\set{x_1,\ldots,x_4}\subset\Z^3$ be the lattice point
vertices of a tetrahedron containing the origin. Let
$\mu_1,\ldots,\mu_4\in\Q$ give the (unique) \emph{barycentric
coordinate} of the origin with respect to the $x_i$;
\begin{align*}
\text{ie. }&\mu_1 x_1+\ldots+\mu_4 x_4=0\\
&\mu_1+\ldots+\mu_4=1\\
&\mu_1\geq 0,\ldots,\mu_4\geq 0.
\end{align*}
Choose $\lambda_1,\ldots,\lambda_4\in\N$ coprime such that
$\mu_i=\frac{\lambda_i}{h}$, where $h=\lambda_1+\ldots+\lambda_4$.
\begin{defn}
Let $q\in\Q$. We define $\intpart{q}:=\max\set{a\in\Z\mid a\leq
q}$ and $\roof{q}:=\min\set{a\in\Z\mid a\geq q}$. The
\emph{fractional part} of $q$, denoted $\fracpart{q}$, is given by
$q-\intpart{q}$.
\end{defn}
\begin{lem}\label{narrowLem}
For any $\kappa\in\set{2,\ldots,h-2}$ we have that
$\sum_{i=1}^4\modpart{i}\in\set{1,2,3}$.
\end{lem}
\begin{proof}
Since $\sum_{i=1}^4\lfrac{i}=\kappa\in\N$ it follows that
$\sum_{i=1}^4\modpart{i}\in\set{0,1,2,3}$. Suppose for some
$\kappa\in\set{2,\ldots,h-2}$, $\modpart{i}=0$ for $i=1,2,3,4$. We
have that $h\mid\kappa\lambda_i$ for each $i$, so let $p$ be a
prime such that $p\mid h$, so that $h=p^rh'$ where $p\nmid h'$.
Then $p^r\mid\kappa\lambda_i$. Suppose that $p^r\nmid\kappa$. Then
$p\mid\lambda_i$ for each $i$. Hence $p\mid
(\lambda_1,\ldots,\lambda_4)=1$, a contradiction. Thus
$p^r\mid\kappa$. By induction on the prime divisors of $h$ we see
that $h\mid\kappa$, so in particular $h\leq\kappa$, which is a
contradiction.
\end{proof}
For convenience we make the following definition:
\begin{defn}\label{tetAdmiss}
We say a tetrahedron is \emph{Fano} if the vertices lie at lattice
points and the only non-vertex lattice point it contains is the
origin, which lies strictly in the interior of the tetrahedron.
\end{defn}
\begin{prop}\label{Tet2Sum}
If the tetrahedron associated with the $\lambda_i$ is Fano then
\begin{enumerate}
 \item[(i)]$\sum_{i=1}^4\modpart{i}=2$ for all
$\kappa\in\set{2,\ldots,h-2}$
 \item[and (ii)]$\coprime{\lambda_i,\lambda_j}=1$ for $i\neq j$.
\end{enumerate}
\end{prop}
\begin{proof}
Let the $\lambda_i$ be associated with a Fano tetrahedron. Since
the origin is strictly in the interior the $\lambda_i$ are all
non-zero. By Lemma~\ref{narrowLem} we only need to consider the
cases when $\sum_{i=1}^4\modpart{i}=1$ or $3$. But if
$\sum_{i=1}^4\modpart{i}=3$ for some $\kappa$, then
$\sum_{i=1}^4\invmodpart{i}=1$.

Suppose for some $\kappa\in\set{2,\ldots,h-2}$ the sum is $1$. Let
$\chi_i=\modpart{i}$. Then $(\chi_1,\ldots,\chi_4)$ is the
(unique) barycentric coordinate for some point in the tetrahedron.
We shall show that it is a non-vertex lattice point not equal to
the origin.

We have that $\sum_{i=1}^4\wholepart{i}x_i$ is a lattice point,
call it $a\in\Z^3$. We also have that
$\sum_{i=1}^4\lfrac{i}x_i=0$. Thus
$$\sum_{i=1}^4\chi_i
x_i=\sum_{i=1}^4\lfrac{i}x_i-\sum_{i=1}^4\wholepart{i}x_i=-a\in\Z^3
$$
By the uniqueness of barycentric coordinates we have that $-a$ is
a non-vertex point, since each $\chi_i<1$. Furthermore suppose
$-a=0$, so that $\chi_i=\lambda_i$ for $i=1,2,3,4$. For each $i$,
$\lfrac{i}-\wholepart{i}=\frac{\lambda_i}{h}$, so we obtain that
$\wholepart{i}=\frac{\lambda_i(\kappa-1)}{h}$ and hence that
$h\mid\lambda_i(\kappa-1)$. As in the proof of
Lemma~\ref{narrowLem} we find that $h\mid\kappa-1$, and so in
particular $h+1\leq\kappa$. This contradicts our range for
$\kappa$. Hence $-a$ must be a non-vertex, non-zero lattice point
in the tetrahedron, contradicting our hypothesis.

Now suppose for a contradiction that
$\coprime{\lambda_1,\lambda_2}\neq 1$. We have
$$\frac{\lambda_3}{\coprime{\lambda_1,\lambda_2}}x_3+\frac{\lambda_4}{\coprime{\lambda_1,\lambda_2}}x_4=-\frac{\lambda_1}{\coprime{\lambda_1,\lambda_2}}x_1-\frac{\lambda_2}{\coprime{\lambda_1,\lambda_2}}x_2\in\Z^3.$$
Since the triangle with vertices $\set{0,x_3,x_4}$ is lattice
point free, by Lemma~\ref{emptyTriEqu} there exists an element of
$GL(3,\Z)$ mapping $x_3\mapsto e_1$ and $x_4\mapsto e_2$. Hence it
must be that $\coprime{\lambda_1,\lambda_2}\mid\lambda_3$ and
$\coprime{\lambda_1,\lambda_2}\mid\lambda_4$, thus
$\coprime{\lambda_1,\ldots,\lambda_4}\neq 1$.
\end{proof}
\begin{cor}\label{ResH}
Let $\coord{\lambda_1,\ldots,\lambda_4}$ be associated with a Fano
tetrahedron. Then
\begin{enumerate}
 \item[(i)]$\sum_{i=1}^4\rfrac{i}=\kappa+2$ for all $\kappa\in\set{2,\ldots,h-2}$
 \item[and (ii)]$\coprime{\lambda_i,h}=1$ for $i=1,\ldots,4$.
\end{enumerate}
\end{cor}
\begin{proof}
Proposition~\ref{Tet2Sum} tells us that for a Fano tetrahedron
\begin{align*}
\modsum{1}+\ldots+\modsum{4}&=2h\\
&=2(\lambda_1+\ldots+\lambda_4).
\end{align*}
Thus for each $\kappa\in\set{2,\ldots,h-2}$ and each
$i=1,\ldots,4$ there exists $\con{i}\in\N$ satisfying
\begin{equation}\label{eq1}
0<\kappa\lambda_i-\con{i}h\leq h
\end{equation}
\begin{equation}\label{eq2}
\con{1}+\ldots+\con{4}=\kappa-2.
\end{equation}
Rearranging (\ref{eq1}) we get an expression for the $\con{i}$,
namely
$\frac{\kappa}{h}\lambda_i-1\leq\con{i}<\frac{\kappa}{h}\lambda_i$,
and taking into account the fact that the $\con{i}$ are integers
gives us $\con{i}=\rfrac{i}-1$. Combining this result with
(\ref{eq2}) gives
\begin{equation}\label{eq3}
\rfrac{1}+\ldots+\rfrac{4}=\kappa+2.
\end{equation}
By hypothesis we have $\modpart{1}+\ldots+\modpart{4}=2$, and by
definition
$\frac{\lambda_1\kappa}{h}+\ldots+\frac{\lambda_4\kappa}{h}=\kappa$.
Hence we obtain $\wholepart{1}+\ldots+\wholepart{4}=\kappa-2$.
This, combined with (\ref{eq3}), tells us that $\modpart{i}>0$.
Finally suppose that, for some $i$, $\coprime{\lambda_i,h}\neq 1$.
Then taking
$\kappa=\frac{h}{\coprime{\lambda_i,h}}\in\set{2,\ldots,h-2}$ we
have $\modpart{i}=0$. Hence $\coprime{\lambda_i,h}=1$.
\end{proof}
Although not required, it is worth observing the similarity
between Corollary~\ref{ResH} and the following:
\begin{prop}[\cite{Sca85}]
Let a lattice point tetrahedron containing no non-vertex lattice
points have the vertices of Lemma~\ref{basisTet} with $x,y,z\geq
1$.  Let $d:=x+y+z-1$. Then
\begin{enumerate}
 \item[(i)]$\roof{\frac{\kappa x}{d}}+\roof{\frac{\kappa y}{d}}+\roof{\frac{\kappa z}{d}} = \kappa+2$ for all $\kappa\in\set{1,\ldots,d-1}$
 \item[and (ii)]$\coprime{x,d}=\coprime{y,d}=\coprime{z,d}=1$.
\end{enumerate}
\end{prop}
Let $h\geq 4$. By making use of Corollary~\ref{ResH} we can
construct bounds on the $\lambda_i$. We may assume without loss of
generality that $\lambda_1\leq\ldots\leq\lambda_4$. For each
$\kappa\in\set{2,\ldots,h-2}$ and each $i$ let
$\ai{i}{\kappa}=\rfrac{i}$. The following conditions are immediate
\begin{equation}\label{reCon}
\begin{array}{l}
\ai{1}{\kappa}\leq\ldots\leq\ai{4}{\kappa}\\
\ai{1}{\kappa}+\ldots+\ai{4}{\kappa}=\kappa+2\\
(\ai{1}{2},\ai{2}{2},\ai{3}{2},\ai{4}{2})=\coord{1,1,1,1}.\\
\end{array}
\end{equation}
We have also that
$\frac{h}{\kappa}(\ai{i}{\kappa}-1)<\lambda_i<\frac{h}{\kappa}\ai{i}{\kappa}$,
and so
$$h\max_{2\leq
n\leq\kappa}{\frac{1}{n}(\ai{i}{n}-1)}<\lambda_i<h\min_{2\leq
n\leq\kappa}{\frac{1}{n}\ai{i}{n}}.$$ Recalling
that$\frac{\lambda_i}{h}=\mu_i$ gives us
\begin{equation}\label{muIneq}
\frac{1}{\kappa}(\ai{i}{\kappa}-1)<\mu_i<\frac{1}{\kappa}\ai{i}{\kappa}\\
\end{equation}
\begin{equation}\label{reIneq}
\max_{2\leq
n\leq\kappa}{\frac{1}{n}(\ai{i}{n}-1)}<\mu_i<\min_{2\leq
n\leq\kappa}{\frac{1}{n}\ai{i}{n}}.
\end{equation}

This suggests a recursive method of determining an upper bound on
$h$. Assume $h\geq 4$ is associated with a Fano tetrahedron. Then
it is possible to construct a sequence
$\set{(\ai{1}{\kappa},\ldots,\ai{4}{\kappa})}_{2\leq \kappa\leq
h-2}$ satisfying conditions (\ref{reCon}) and (\ref{reIneq}) for
all $\kappa\in\set{2,\ldots,h-2}$. Moreover we have that for each
$\kappa\in\set{2,\ldots,h-3}$ there exists some
$i\in\set{1,\ldots,4}$ such that
$$\ai{j}{\kappa+1}=\left\{
\begin{array}{ll}
\ai{j}{\kappa}&\text{for }j\neq i\\
\ai{j}{\kappa}+1&\text{for }j=i.
\end{array}\right.$$
\begin{lem}\label{ineqLem}
Let $a,k\in\N$ be such that $a<k$. Then
$\frac{a}{k}>\frac{a}{k+1}$ and $\frac{a}{k}<\frac{a+1}{k+1}$.
\end{lem}
An immediate consequence of Lemma~\ref{ineqLem} is that
$$\frac{1}{\kappa+1}\ai{i}{\kappa+1}=\frac{1}{\kappa+1}(\ai{i}{\kappa}+1)\geq\min_{2\leq n\leq\kappa}{\frac{1}{n}\ai{i}{n}}$$
and hence, using (\ref{muIneq}) and (\ref{reIneq}), we obtain
$$\frac{1}{\kappa+1}(\ai{i}{\kappa+1}-1)=\frac{1}{\kappa+1}\ai{i}{\kappa}<\mu_i<\min_{2\leq
n\leq\kappa+1}{\frac{1}{n}\ai{i}{n}}=\min_{2\leq
n\leq\kappa}{\frac{1}{n}\ai{i}{n}}.$$ Thus we have the requirement
that
\begin{equation}\label{testIneq}
\frac{1}{\kappa+1}\ai{i}{\kappa}<\min_{2\leq
n\leq\kappa}{\frac{1}{n}\ai{i}{n}}.
\end{equation}

Conditions (\ref{reCon}) and (\ref{testIneq}) are independent of
$h$, so by writing a simple recursive function on a computer it is
possible to test these conditions for large values of $\kappa$,
using all the sequences obtained for $\kappa$ to check whether a
sequence exists for $\kappa+1$. If no such sequence exists we have
found an upper bound on $h$, namely $h\leq\kappa+2$.

It is worth observing that this method for finding a bound on $h$
really does do that; when all possible sequences have terminated
it is impossible to proceed any further. No \emph{a priori}
guarantee that this search along all possible sequences will
terminate has been given here.

It is also worth noting that the bound this method gives is not
the tightest, however this deficiency is balanced by the fact that
it providing a technique which is independent of $h$.

This yields a bound of $h\leq 30$. Proposition~\ref{mainProp} now
follows from Proposition~\ref{Tet2Sum} by the easy task of
checking all possible $\lambda_i$ up to this bound. An alternative
proof of Proposition~\ref{mainProp} can be found in~\cite{BB92}.
\begin{prop}\label{mainProp}
Let $\lambda_1\leq\ldots\leq\lambda_4$ be associated with a Fano
tetrahedron. Then $\coord{\lambda_1,\ldots,\lambda_4}$ is equal to
one of the following:
\begin{center}
\begin{tabular}{cccc}
$\coord{1,1,1,1}$&$\coord{1,1,1,2}$&$\coord{1,1,2,3}$&$\coord{1,2,3,5}$\\
$\coord{1,3,4,5}$&$\coord{2,3,5,7}$&$\coord{3,4,5,7}$.&
\end{tabular}
\end{center}
\end{prop}
\section{Classifying the Tetrahedra: The Coordinates of the Vertices}
\begin{prop}\label{tetCoordRel}
Let $\lambda_1\leq\ldots\leq\lambda_4$ be associated with a Fano
tetrahedron. Then, by means of the action of $GL(3,\Z)$, we can
transform the tetrahedron to the form
$$\begin{pmatrix}1&0&k''\lambda_4-a\lambda_1&-k''\lambda_3-b\lambda_1\\
0&1&k'\lambda_4-a\lambda_2&-k'\lambda_3-b\lambda_2\\
0&0&k\lambda_4&-k\lambda_3\\\end{pmatrix}$$ where $a,b\in\Z, a>0$
are such that $a\lambda_3+b\lambda_4=1$, and $k,k',k''\in\N$ are
such that
\begin{align*}&0\leq k''\lambda_4-a\lambda_1<k\lambda_4\\
\text{and }&0\leq k'\lambda_4-a\lambda_2<k\lambda_4
\end{align*}
with one of these inequalities equal to zero only if
$\lambda_4=1$.
\end{prop}
\begin{proof}
By virtue of Lemma~\ref{emptyTriEqu} we may assume without loss of
generality that our tetrahedron has vertices $\set{e_1,e_2,x,y}$
with $\lambda_1e_1+\lambda_2e_2+\lambda_3x+\lambda_4y=0$. Thus we
see that $\lambda_3x_3=-\lambda_4y_3$, and so
$y_3=-\frac{\lambda_3}{\lambda_4}x_3\in\Z$. Hence
$\lambda_4\mid\lambda_3x_3$, but $\coprime{\lambda_3,\lambda_4}=1$
and so it must be that $\lambda_4\mid x_3$. Thus there exists some
$k\in\Z$ such that
\begin{align*}
x_3&=k\lambda_4\\
y_3&=-k\lambda_3.
\end{align*}
We may take $x_3$ positive, and so $k\in\N$.

We also have that $\lambda_2+\lambda_3x_2+\lambda_4y_2=0$, so that
$\lambda_3x_2+\lambda_4y_2=-\lambda_2$. Now since
$\coprime{\lambda_3,\lambda_4}=1$ there exists $a,b\in\Z, a>0$
such that $\lambda_3a+\lambda_4b=1$. This gives us that
$\lambda_3(-\lambda_2a)+\lambda_4(-\lambda_2b)=-\lambda_2$, so
that $\lambda_3(x_2+a\lambda_2)+\lambda_4(y_2+b\lambda_2)=0$. Thus
there exists some $k'\in\Z$ such that
\begin{align*}
x_2&=k'\lambda_4-a\lambda_2\\
y_2&=-k'\lambda_3-b\lambda_2.
\end{align*}
Similarly we obtain that there exists some $k''\in\Z$ such that
\begin{align*}
x_1&=k''\lambda_4-a\lambda_1\\
y_1&=-k''\lambda_3-b\lambda_1.
\end{align*}

By applying
$$\begin{pmatrix}1&0&c\\0&1&d\\0&0&1\end{pmatrix}\in GL(3\Z)$$
for suitably chosen $c,d\in\Z$ we can arrange matters such that
(with possible relabeling of $k'$ and $k''$)
\begin{align*}
0&\leq k'\lambda_4-a\lambda_2<k\lambda_4\\
0&\leq k''\lambda_4-a\lambda_1<k\lambda_4.
\end{align*}

Now suppose that $k'\lambda_4-a\lambda_2=0$. Since
$\coprime{\lambda_2,\lambda_4}=1$ there must exist some constant
$m\in\Z$ such that $k'=m\lambda_2$ and $a=m\lambda_4$. In
particular this gives us that $\lambda_4(m\lambda_3+b)=1$, so that
$\lambda_4=1$. Similarly if $k''\lambda_4-a\lambda_1=0$. This
completes the proof.
\end{proof}
The exceptional case in Proposition~\ref{tetCoordRel} occurring
when $\lambda_1=\ldots=\lambda_4=1$ will be dealt with now.
\begin{prop}\label{exceptionalTet}
With notation as above, the only exceptional case is given, up to
equivalence, by the tetrahedron with vertices
$\set{e_1,e_2,e_3,-e_1-e_2-e_3}$.
\end{prop}
\begin{proof}
With notation as in the proof of Proposition~\ref{tetCoordRel}, we
have that $a=1,b=0$ and so taking $k'\lambda_4-a\lambda_2=0$
implies that $k'=1$. Thus we see that our tetrahedron has the form
$$\begin{pmatrix}1&0&k''-1&-k''\\0&1&0&-1\\0&0&k&-k\\\end{pmatrix}$$
where $k''$ and $k$ are to be determined.

The triangle defined by the origin, the first and the third
vertices in the above matrix is lattice point free. Thus
$$\det\begin{pmatrix}1&k''-1\\0&k\\\end{pmatrix}=\pm1.$$
This forces $k=1$ and the resulting tetrahedron is equivalent to
that given in the statement.
\end{proof}
The following two results are taken from \cite{Sca85}. A proof is
given for the first result because we need to know explicitly the
steps required for the transformation.
\begin{lem}[\cite{Sca85}]\label{basisTet}
A lattice point tetrahedron containing no non-vertex lattice
points can, by means of a translation and the action of
$GL(3,\Z)$, be transformed to the form
$$\begin{pmatrix}1&0&0&x\\0&1&0&y\\0&0&1&z\\\end{pmatrix}$$
where $x,y,z\in\Z, x,y\geq 0, z\geq 1$.
\end{lem}
\begin{proof}
By applying a translation if necessary and considering
Lemma~\ref{emptyTriEqu}, we may assume without loss of generality
that the tetrahedron is in the form
$$\begin{pmatrix}0&1&0&x\\0&0&1&y\\0&0&0&z\\\end{pmatrix}$$
where $z\geq 1$, but conditions on $x$ and $y$ still to be
determined.

Let $x\mapsto x ~(\text{mod }{z})$ and $y\mapsto y ~(\text{mod
}{z})$. Observe that this is equivalent to the (left) action of
$$\begin{pmatrix}1&0&a\\0&1&b\\0&0&1\end{pmatrix}\in GL(3,\Z)$$
for suitably chosen $a,b\in\Z$. Thus we can assume that $0\leq
x<z$ and $0\leq y<z$. Suppose that $z<x+y$. Then  set
$$\begin{array}{llll}
\mu_1:=1-\mu_2-\mu_3-\mu_4,&\mu_2:=1-\frac{x}{z},&\mu_3:=1-\frac{y}{z},&\mu_4:=\frac{1}{z}.\end{array}$$
Clearly $\sum\mu_i=1$, and $\mu_2,\mu_3,\mu_4\geq 0$. We have also
that $\mu_1=\frac{1}{z}(x+y-z-1)\geq 0$. But then
$$\mu_1\begin{pmatrix}0\\0\\0\\\end{pmatrix}+\mu_2\begin{pmatrix}1\\0\\0\\\end{pmatrix}+\mu_3\begin{pmatrix}0\\1\\0\\\end{pmatrix}+\mu_4\begin{pmatrix}x\\y\\z\\\end{pmatrix}=\begin{pmatrix}1\\1\\1\\\end{pmatrix}$$
gives us a non-vertex lattice point in the interior of the
tetrahedron - a contradiction. Thus it must be that $z\geq x+y$.

Finally we apply the unimodular transformation $z\mapsto -x-y+z+1$
which yields the result.
\end{proof}
\begin{prop}[\cite{Sca85} Simplification of Howe's Theorem]\label{Howe}
Let a lattice point tetrahedron containing no non-vertex lattice
points have the vertices of Lemma~\ref{basisTet} with $x,y,z\geq
1$. Then $\set{x,y,z}\cap\set{1}\neq\emptyset$.
\end{prop}
Let us now consider a Fano tetrahedron presented in the form given
in Proposition~\ref{tetCoordRel}. In addition we shall assume we
are not looking at the case handled in
Proposition~\ref{exceptionalTet}. The tetrahedron with vertices
given by $\set{0,e_1,e_2,(x,y,z)}$, where
$$\begin{array}{lll}
x:=k''\lambda_4-a\lambda_1\geq 1,& y:=k'\lambda_4-a\lambda_2\geq
1,& z:=k\lambda_4\geq 1,
\end{array}$$
is lattice point free. By following the proof of
Lemma~\ref{basisTet} we see that it is equivalent to
$$\begin{pmatrix}1&0&0&x\\0&1&0&y\\0&0&1&z-x-y+1\\\end{pmatrix}$$
and that $z\geq x+y$. Proposition~\ref{Howe} tells us that
$\set{x,y,z-x-y+1}\cap\set{1}\neq\emptyset$. Thus
$$\begin{array}{rcl}
\text{either }k''=\frac{1+a\lambda_1}{\lambda_4}\in\Z&\text{if and only if}&x=1\\
\text{or }k'=\frac{1+a\lambda_2}{\lambda_4}\in\Z&\text{if and only if}&y=1\\
\text{or
}k-k'-k''=-a\frac{\lambda_1+\lambda_2}{\lambda_4}\in\Z&\text{if
and only if}&z-x-y+1=1.
\end{array}$$
The result of applying this to the barycentric coordinates found
in Proposition~\ref{mainProp} is given in Table~\ref{partVertTab}.
Observe that the only cases of ambiguity are for $\coord{1,1,1,1}$
and $\coord{1,1,1,2}$.
\begin{table}[ht]
\centering
\begin{tabular}{|c|c|c|c|c|c|}
\hline
$\coord{\lambda_1,\lambda_2,\lambda_3,\lambda_4}$&$a$&$b$&$\frac{1+a\lambda_1}{\lambda_4}$&$\frac{1+a\lambda_2}{\lambda_4}$&$a\frac{\lambda_1+\lambda_2}{\lambda_4}$\\
\hline\hline
$\coord{1,1,1,1}$&1&0&\multicolumn{2}{|c|}{2}&2\\
\hline
$\coord{1,1,1,2}$&1&0&\multicolumn{2}{|c|}{1}&1\\
\hline
$\coord{1,1,2,3}$&2&-1&\multicolumn{2}{|c|}{1}&-\\
\hline
$\coord{1,2,3,5}$&2&-1&-&1&-\\
\hline
$\coord{1,3,4,5}$&4&-3&1&-&-\\
\hline
$\coord{2,3,5,7}$&3&-2&1&-&-\\
\hline
$\coord{3,4,5,7}$&3&-2&-&-&3\\
\hline
\end{tabular}
\caption{The values depending on $a$.} \label{partVertTab}
\end{table}
\begin{prop}\label{otherSizeTet}
Let $\lambda_1\leq\ldots\leq\lambda_4$ be associated with a Fano
tetrahedron presented in the form given in
Proposition~\ref{tetCoordRel}. Then
\begin{align*}
&0\leq k\lambda_3-k''\lambda_3-b\lambda_1<k\lambda_3\\
\text{and }&0\leq k\lambda_3-k'\lambda_3-b\lambda_2<k\lambda_3
\end{align*}
with one of these inequalities equal to zero only if
$\lambda_3=1$, in which case the tetrahedron is equivalent either
to that given in Proposition~\ref{exceptionalTet} or to
$$\begin{pmatrix}1&0&1&-1\\0&1&1&-1\\0&0&2&-1\\\end{pmatrix}.$$
\end{prop}
\begin{proof}
Since $a\lambda_3+b\lambda_4=1$ we have that
$a=\frac{1-b\lambda_4}{\lambda_3}$. Substituting this into the
inequalities obtained in Proposition~\ref{tetCoordRel} gives us
that $0\leq k\lambda_3-k''\lambda_3-b\lambda_1<k\lambda_3$ and
$0\leq k\lambda_3-k'\lambda_3-b\lambda_2<k\lambda_3$. Now suppose
that $k\lambda_3-k''\lambda_3-b\lambda_1=0$. Then we have that
$(k-k'')\lambda_3=b\lambda_1$, and since
$\coprime{\lambda_1,\lambda_3}=1$ there must exist some $c\in\Z$
such that $k-k''=c\lambda_1$ and $b=c\lambda_3$. But then
$a\lambda_3+c\lambda_3\lambda_4=1$ which forces $\lambda_3=1$ (as
required). The only cases when $\lambda_3=1$ are when $a=1,b=0$.
Hence $k=k''$.

There are two possible choices for $\lambda_4$. First consider the
case when $\lambda_4=1$. We have that $k\geq k''+k'-2$, and
$k'\geq 2$. Thus $k'=2$. Hence we see that our Fano tetrahedron is
equivalent to the form
$$\begin{pmatrix}1&0&-1&0\\0&1&1&-2\\0&0&k&-k\\\end{pmatrix}.$$
The triangle with vertices given by the origin, and the second and
fourth column of the above matrix is lattice point free. By
Lemma~\ref{emptyTriEqu} it must be that $k=1$, which gives a
tetrahedron equivalent to that derived in
Proposition~\ref{exceptionalTet}.

Now consider the case when $\lambda_4=2$. We have that $k\geq
k''+k'-1$, and $k'\geq 1$. Thus $k'=1$. Hence we see that our Fano
tetrahedron equivalent to the form
$$\begin{pmatrix}1&0&-1&0\\0&1&1&-1\\0&0&2k&-k\\\end{pmatrix}.$$
As before we see that $k=1$ and the result follows.
\end{proof}
We consider a Fano tetrahedron presented in the form given in
Proposition~\ref{tetCoordRel} and assume we are not looking at the
case handled in Proposition~\ref{otherSizeTet}. By
Proposition~\ref{Howe} we have that
\begin{align*}
\text{either }k-k''&=\frac{1+b\lambda_1}{\lambda_3}\in\Z\\
\text{or }k-k'&=\frac{1+b\lambda_2}{\lambda_3}\in\Z\\
\text{or }k-k'-k''&=b\frac{\lambda_1+\lambda_2}{\lambda_3}\in\Z.
\end{align*}
The result of applying this to the barycentric coordinates found
in Proposition~\ref{mainProp} is presented in
Table~\ref{part2VertTab}. The results of Table~\ref{partVertTab}
and Table~\ref{part2VertTab} complement each other beautifully,
allowing the relationships amongst $k,k'$ and $k''$ shown in
Table~\ref{vertRelTab} to be established.
\begin{table}[ht]
\centering
\begin{tabular}{|c|c|c|c|c|c|}
\hline
$\coord{\lambda_1,\lambda_2,\lambda_3,\lambda_4}$&$a$&$b$&$\frac{1+b\lambda_1}{\lambda_3}$&$\frac{1+b\lambda_2}{\lambda_3}$&$b\frac{\lambda_1+\lambda_2}{\lambda_3}$\\
\hline\hline
$\coord{1,1,1,1}$&1&0&\multicolumn{2}{|c|}{1}&0\\
\hline
$\coord{1,1,1,2}$&1&0&\multicolumn{2}{|c|}{1}&0\\
\hline
$\coord{1,1,2,3}$&2&-1&\multicolumn{2}{|c|}{0}&-1\\
\hline
$\coord{1,2,3,5}$&2&-1&0&-&-1\\
\hline
$\coord{1,3,4,5}$&4&-3&-&-2&-3\\
\hline
$\coord{2,3,5,7}$&3&-2&-&-1&-2\\
\hline
$\coord{3,4,5,7}$&3&-2&-1&-&-\\
\hline
\end{tabular}
\caption{The values depending on $b$.} \label{part2VertTab}
\end{table}
\begin{table}[ht]
\centering
\begin{tabular}{|c|c|c|c|}
\hline
$\coord{\lambda_1,\lambda_2,\lambda_3,\lambda_4}$&$k$&$k'$&$k''$\\
\hline \hline
$\coord{1,1,1,1}$&$k$&$k-2$&2\\
\cline{2-4} &$k$&3&$k-1$\\
\hline
$\coord{1,1,1,2}$&$k$&$k-1$&1\\
\cline{2-4} &$k$&$k-1$&2\\
\hline $\coord{1,1,2,3}$&$k$&$k$&1\\
\cline{2-4} &1&1&1\\
\hline
$\coord{1,2,3,5}$&$k$&1&$k$\\
\hline
$\coord{1,3,4,5}$&$k$&$k-2$&1\\
\cline{2-4} &$k$&$k+2$&1\\
\hline
$\coord{2,3,5,7}$&$k$&$k+1$&1\\
\hline
$\coord{3,4,5,7}$&$k$&2&$k+1$\\
\hline
\end{tabular}
\caption{The relationships between $k,k'$ and $k''$.}
\label{vertRelTab}
\end{table}

We are now in a position to calculate the vertices of the Fano
tetrahedra (up to the action of $GL(3,\Z)$). We will proceed by
taking each barycentric coordinate in turn and combining the
results of Table~\ref{vertRelTab} and
Proposition~\ref{tetCoordRel}. The final results are collected
together in Table~\ref{tetTable}. It is worth comparing this with
the results of~\cite{Suz02}.

For ease of notation in the sequel, we shall make the temporary
convention that by $x_i$ we mean the $i^\text{th}$ column of the
matrix $\begin{pmatrix}x_1^t&x_2^t&x_3^t&x_4^t\end{pmatrix}$
representing a tetrahedron with vertices equivalent to
$\set{x_1,\ldots,x_4}$.

First we consider the case with barycentric coordinate
$\coord{1,1,1,1}$. From the results of Table~\ref{vertRelTab} and
Proposition~\ref{tetCoordRel} we have that our Fano tetrahedron
has two possible forms, both of are equivalent to
$$\begin{pmatrix}1&0&1&-2\\0&1&-3&2\\0&0&k&-k\\\end{pmatrix}.$$
We observe that $x_3$ tells us that $\coprime{3,k}=1$ and $x_4$
tells us that $\coprime{2,k}=1$. Furthermore, taking $k=1$ gives
us a tetrahedron equivalent to that found in
Proposition~\ref{exceptionalTet}. Suppose that $k\geq 7$. Then
$\frac{4}{k}x_2+\frac{2}{k}x_3+\frac{1}{k}x_4=e_3$, which
contradicts our tetrahedron being Fano. Thus the only remaining
possibility is that $k=5$, which by inspection we see does indeed
give us a Fano tetrahedron.

Now we consider the case with barycentric coordinate
$\coord{1,1,1,2}$. By Table~\ref{vertRelTab} and
Proposition~\ref{tetCoordRel} we see once more that our Fano
tetrahedron can take two possible forms. First we consider the
form equivalent to
$$\begin{pmatrix}1&0&1&-1\\0&1&-3&1\\0&0&2k&-k\\\end{pmatrix}.$$
If we take $k=1$ we obtain a Fano tetrahedron equivalent to that
found in Proposition~\ref{otherSizeTet}. Suppose that $k=2$. Then
$\frac{1}{2}\coord{1,-3,4}+\frac{1}{2}\coord{-1,1,-2}=\coord{0,-1,1}$
is a non-vertex, non-zero lattice point in the interior of the
tetrahedron, and hence it is not Fano. The third column tells us
that $\coprime{3,k}=1$. Finally the tetrahedron is not Fano if
$k\geq 4$ since
$\frac{2}{k}x_2+\frac{1}{k}x_3+\frac{1}{k}x_4=e_3$.

Now we consider the second possibility, which is equivalent to
$$\begin{pmatrix}1&0&3&-2\\0&1&-3&1\\0&0&2k&-k\\\end{pmatrix}.$$
When $k=1$ we obtain a Fano tetrahedron equivalent to the one
previously found. $x_3$ and $x_4$ tell us that $\coprime{3,k}=1$
and $\coprime{2,k}=1$ respectively, and if $k\geq 7$ we have the
non-vertex, non-zero internal lattice point given by
$\frac{3}{k}x_1+\frac{1}{k}x_3+\frac{3}{k}x_4=-e_3$. Thus the only
possibility left is $k=5$, which contains the lattice point
$\frac{1}{5}\coord{1,0,0}+\frac{2}{5}\coord{3,-3,10}+\frac{1}{5}\coord{-2,1,-5}=\coord{1,-1,3}$.

For barycentric coordinate $\coord{1,1,2,3}$ the two possibilities
are (up to equivalence)
$$\begin{pmatrix}1&0&1&-1\\0&1&-2&1\\0&0&3k&-2k\\\end{pmatrix}\text{ and }\begin{pmatrix}1&0&1&-1\\0&1&1&-1\\0&0&3&-2\\\end{pmatrix}.$$
The third column tells us that $k$ must be odd, but if $k\geq 3$
we have the interior lattice point
$\frac{1}{k}x_2+\frac{1}{k}x_3+\frac{1}{k}x_4=e_3$. Thus the only
possibility is that $k=1$, but the resulting tetrahedron is
equivalent to that already found.

When we have barycentric coordinate $\coord{1,2,3,5}$ our
tetrahedron equivalent to
$$\begin{pmatrix}1&0&-2&1\\0&1&1&-1\\0&0&5k&-3k\\\end{pmatrix}.$$
The third column tells us that $k$ is odd, and if $k\geq 3$ we
have the internal lattice point
$\frac{1}{k}x_1+\frac{1}{k}x_3+\frac{1}{k}x_4=2e_3$. By inspection
we see that the case when $k=1$ is Fano.

For barycentric coordinate $\coord{1,3,4,5}$ we have two
possibilities. First we consider the case when our tetrahedron is
equivalent to
$$\begin{pmatrix}1&0&1&-1\\0&1&-22&17\\0&0&5k&-4k\\\end{pmatrix}.$$
$x_3$ tells us that $k$ is odd. If $k\geq 7$ then it is not Fano,
since $\frac{5}{k}x_2+\frac{1}{k}x_3+\frac{1}{k}x_4=e_3$. If $k=5$
then
$\frac{1}{5}\coord{1,-22,25}+\frac{1}{5}\coord{-1,17,-20}=\coord{0,-1,1}$,
and if $k=3$ then $k=3$ then
$\frac{1}{3}e_1+\frac{1}{3}e_2+\frac{1}{3}\coord{-1,17,-12}=\coord{0,6,-4}$.
By inspection we see that the case when $k=1$ is Fano.

The second possibility is the tetrahedron equivalent to
$$\begin{pmatrix}1&0&1&-1\\0&1&-2&1\\0&0&5k&-4k\\\end{pmatrix}.$$
We require that $k$ is odd, but if $k\geq 3$ we obtain the point
$\frac{1}{k}x_2+\frac{1}{k}x_3+\frac{1}{k}x_4=e_1$, and when $k=1$
we obtain the tetrahedron found above.

Continuing in the same vein, for barycentric coordinate
$\coord{2,3,5,7}$ we have
$$\begin{pmatrix}1&0&1&-1\\0&1&-2&1\\0&0&7k&-5k\\\end{pmatrix}.$$
This tells us that $k$ is odd, and if $k\geq 3$ we obtain the
internal lattice point
$\frac{1}{k}x_2+\frac{1}{k}x_3+\frac{1}{k}x_4=2e_3$. Thus $k=1$ is
the only possibility, and we see by inspection that it is indeed
Fano.

Finally consider barycentric coordinate $\coord{3,4,5,7}$. This
gives us
$$\begin{pmatrix}1&0&-2&1\\0&1&2&-2\\0&0&7k&-5k\\\end{pmatrix}.$$
Once more we see that $k$ must be odd, and that if $k\geq 3$ then
it is not Fano sice we have
$\frac{1}{k}x_1+\frac{1}{k}x_3+\frac{1}{k}x_4=2e_3$. When $k=1$ we
do indeed get a Fano tetrahedron.
\begin{table}[ht]
\centering
\begin{tabular}{|c|c|c|c|}
\hline
$\frac{1}{4}\coord{1,1,1,1}$&$\frac{1}{4}\coord{1,1,1,1}$&$\frac{1}{5}\coord{1,1,1,2}$&$\frac{1}{7}\coord{1,1,2,3}$\\
\hline
$\begin{pmatrix}1&0&0&-1\\0&1&0&-1\\0&0&1&-1\\\end{pmatrix}$&
$\begin{pmatrix}1&0&1&-2\\0&1&-3&2\\0&0&5&-5\\\end{pmatrix}$&
$\begin{pmatrix}1&0&1&-1\\0&1&1&-1\\0&0&2&-1\\\end{pmatrix}$&
$\begin{pmatrix}1&0&1&-1\\0&1&-2&1\\0&0&3&-2\\\end{pmatrix}$\\
\hline\hline
$\frac{1}{11}\coord{1,2,3,5}$&$\frac{1}{13}\coord{1,3,4,5}$&$\frac{1}{17}\coord{2,3,5,7}$&$\frac{1}{19}\coord{3,4,5,7}$\\
\hline
$\begin{pmatrix}1&0&-2&1\\0&1&1&-1\\0&0&5&-3\\\end{pmatrix}$&
$\begin{pmatrix}1&0&1&-1\\0&1&-2&1\\0&0&5&-4\\\end{pmatrix}$&
$\begin{pmatrix}1&0&1&-1\\0&1&-2&1\\0&0&7&-5\\\end{pmatrix}$&
$\begin{pmatrix}1&0&-2&1\\0&1&2&-2\\0&0&7&-5\\\end{pmatrix}$\\
\hline
\end{tabular}
\caption{The vertices of the Fano tetrahedra,up to the action of
$GL(3,\Z)$}\label{tetTable}
\end{table}
\section{Classifying the Minimal Polytopes}\label{classifyMinNow}
We extend Definition \ref{tetAdmiss} to any polytope $P$.
\begin{defn}
We say a lattice point polytope $P$ in $\Z^3$ is \emph{Fano} if
$P$ is convex and the only non-vertex lattice point it contains is
the origin, which lies strictly in the interior of the polytope.
\end{defn}
Given any Fano polytope $P$ with vertices $\set{x_1,\ldots,x_k}$
we make the following definition:
\begin{defn}\label{minimalDef}
We say $P$ is \emph{minimal} if, for all $j\in\set{1,\ldots,k}$,
the polytope $P'$ given by the vertices
$\set{x_1,\ldots,x_k}\setminus\set{x_j}$ is not Fano.
\end{defn}
Let us consider a minimal Fano polytope $P$. Since $0\in P$ there
exist non-coplanar vertices $x_1,\dots,x_4$ of $P$ such that
$0\in\sconv{x_1,\dots,x_4}=:P'$.

Either $P$ is equivalent to one of the tetrahedra in
Table~\ref{tetTable}, or it is not. If it is not, then minimality
gives us that it does not contain a Fano tetrahedron; in
particular $P'$ is not Fano. We assume that this is the case.

Since $P'$ is not a Fano tetrahedron it must be that either the
origin lies on a face of $P'$ or on an edge of $P'$. If the origin
lies on a face of $P'$ then $P$ contains a Fano triangle. Thus
there exist three vertices of $P$ which lie in a plane containing
the origin, and the origin lies strictly in the interior of the
triangle defined by these three points. This possibility will be
discussed in further detail below.

Assume now that $P$ does not contain a Fano triangle. Then it must
be that the origin lies on one of the edges of $P'$, say on the
edge defined by $x_1$ and $x_2$. Since the origin lies in the
strict interior of $P$ there must exist distinct vertices
$y_1,\ldots,y_4$ of $P$ not equal to $x_1$ or $x_2$ such that
$\sconv{x_1,x_2,y_1,y_2}$ is a Fano square and
$\sconv{x_1,x_2,y_1,y_2}$ is a Fano
 octahedron. Minimality gives that $P$ is a Fano
 octahedron, and these will be classified in Lemma~\ref{lemCaseII}.

We return now to consider in more detail the case when $P$
contains a Fano triangle, say that defined by $\set{x_1,x_2,x_3}$.
Since the origin lies in the strict interior of $P$ there must
exist vertices $y_1$ and $y_2$ lying on either side of the plane
containing our Fano triangle. Minimality then gives us that $P$ is
precisely the polygon with vertices $\set{x_1,x_2,x_3,y_1,y_2}$.

Now consider the line passing through the origin and $y_1$. This
line crosses the polytope $P$ at points $y_1\in\Z^3$ and $x$ not
necessarily in $\Z^3$. There are three possible locations for $x$:
\begin{enumerate}
 \item{(i)} $x$ is equal to $y_2$. Then $y_2=-y_1$. These will be
 classified in Lemma \ref{lemCaseI1}.
 \item{(ii)} $x$ lies on the edge with endpoints $\set{x_1,y_2}$.
 Then $\set{x_1,y_1,y_2}$ is a Fano triangle. We use the fact that the
 origin has barycentric coordinate
 $\coord{\frac{1}{3},\frac{1}{3},\frac{1}{3}}$ with respect to
 $\set{x_1,y_1,y_2}$. Thus the line passing through $x_1$ and the origin bisects the line with
 endpoints $\set{y_1,y_2}$ at a point $x'$, say. Now the length of
 the line joining $\set{x_1,0}$ is twice the length of the line
 joining $\set{x',0}$. Similarly by considering the Fano triangle
 $\set{x_1,x_2,x_3}$, the line passing through $x_1$ and the origin bisects the line with endpoints
 $\set{x_2,x_3}$ at a point $x''$, say, and we have that the distance from $x_1$ to the
 origin is twice the length of the line joining the origin to
 $x''$. Hence we see that $\set{x_2,x_3,y_1,y_2}$ are coplanar
 and form a parallelogram. These will be classified in Lemma
 \ref{lemCaseI2}.
 \item{(iii)} $x'$ lies strictly in the interior of the triangle
 $\set{x_1,x_2,y_2}$. But then $\set{x_1,x_2,y_1,y_2}$ is a
 Fano tetrahedron, contradicting our assumption.
\end{enumerate}
\begin{lem}\label{lemCaseII}
The vertices of the minimal Fano octahedra (up to the action of
$GL(3,\Z)$) are given by
$$\begin{pmatrix}1&0&-1&0&0&0\\0&1&0&-1&0&0\\0&0&0&0&1&-1\\\end{pmatrix},
\begin{pmatrix}1&0&-1&0&1&-1\\0&1&0&-1&1&-1\\0&0&0&0&2&-2\\\end{pmatrix}.$$
\end{lem}
\begin{proof}
By making use of Lemma~\ref{emptyTriEqu} and recalling that $P$
does not contain a Fano triangle, we may take the vertices of $P$
to be $\set{e_1,-e_1,e_2,-e_2,x_1,x_2}$. We observe that
$x_1=-x_2$, otherwise we would have that $P$ contains a Fano
tetrahedron. So take $x=-x_2=x_1=\coord{a,b,c}$. First we shall
show that, without loss of generality, we may take $a,b,c$ such
that
\begin{equation}\label{ineqII}
0\leq a\leq b\leq\frac{c}{2}
\end{equation}
Trivially we may assume that $0\leq a\leq b$. Suppose that
$b>\frac{c}{2}$. Then $b-c>-\frac{c}{2}$ and so $c-b<\frac{c}{2}$.
This process corresponds to the action of $GL(3,\Z)$ transforming
$$\begin{pmatrix}1&0&-1&0&a&-a\\0&1&0&-1&b&-b\\0&0&0&0&c&-c\\\end{pmatrix}\text{ to }\begin{pmatrix}1&0&-1&0&a&-a\\0&-1&0&1&c-b&-(c-b)\\0&0&0&0&c&-c\\\end{pmatrix}.$$
Hence we may assume that inequality (\ref{ineqII}) holds.

Now consider the point $e_3$. Either $x=e_3$ or $e_3$ lies outside
of $P$. The first possibility gives us the first Fano octahedron.
The second possibility tells us that $e_3$ must lie on the
opposite side to the origin of the plane defined by
$\set{-e_1,-e_2,x}$. This plane intersects the $z$-axis at the
point $\coord{0,0,\frac{c}{a+b+1}}$. This gives us that $c\leq
a+b$. Combining this with (\ref{ineqII}) gives us that $b\leq a$
and so $b=a$. This in turn gives us that $c\leq 2b$ and $2b\leq
c$, and so we obtain $2a=2b=c$. Thus (up to the action of
$GL(3,\Z)$) we have that $a=1,b=1,c=2$, which gives us the second
Fano octahedron.
\end{proof}
\begin{lem}\label{lemCaseI1}
If $P$ is a minimal Fano polytope with vertices
$\set{x_1,x_2,x_3,y_1,-y_1}$ such that $\set{x_1,x_2,x_3}$ are the
vertices of a Fano triangle, then $P$ is equal (up to the action
of $GL(3,\Z)$) to one of
$$\begin{pmatrix}1&0&-1&0&0\\0&1&-1&0&0\\0&0&0&1&-1\end{pmatrix},
\begin{pmatrix}1&0&-1&1&-1\\0&1&-1&2&-2\\0&0&0&3&-3\end{pmatrix}.$$
\end{lem}
\begin{proof}
By making use of Lemma~\ref{emptyTriEqu} we may take the vertices
of $P$ to be $\set{e_1,e_2,-e_1-e_2,x,y}$. If $y\neq -x$ then $P$
would contain a Fano tetrahedron, which contradicts minimality.
Let $x=\coord{a,b,c}$. We claim that, without loss of generality,
we may take $a,b,c$ such that $0<a\leq b\leq c$ and
\begin{equation}\label{ineqI1a}
a+b\leq c
\end{equation}
Clearly we can take $0<a\leq b$ and $c>0$. Suppose that $a+b>c$.
Then we have that $(c-a)+(c-b)<c$. By using the fact that $y=-x$
and applying the transformation
$$\begin{pmatrix}1&0&-c\\0&1&-c\\0&0&-1\\\end{pmatrix}\in GL(3,\Z)$$
we see that we may assume that inequality (\ref{ineqI1a}) holds.

Now consider the point $e_3$. Either $x=e_3$ or $e_3$ lies outside
of $P$. The first case gives us the first Fano polytope in the
statement. The second case tells us that we have $e_3$ lies on the
opposite side to the origin of the plane defined by
$\set{e_1,-e_1-e_2,x}$. This plane intersects the $z$-axis at the
point $\coord{0,0,\frac{c}{2b-a+1}}$, and so
\begin{equation}\label{ineqI1b}
2b-a\geq c.
\end{equation}

Now consider the point $x'=e_2+e_3$. Either $x'=x$, which gives a
Fano polytope equivalent to the one previously found, or $x'$ lies
outside of $P$. If this is the case we have that $x'$ lies on the
opposite side to the origin of the plane defined by
$\set{e_2,-e_1-e_2,x}$. This plane intersects the line passing
through the origin and $e_2+e_3$ at the point $\coord{0,k,k}$
where $k:=\frac{c}{2a-b+c+1}$. Hence
\begin{equation}\label{ineqI1c}
b\leq 2a.
\end{equation}

Now suppose both $e_3$ and $x'$ lie outside $P$. Combining
inequalities (\ref{ineqI1a}) and (\ref{ineqI1b}) gives us that
$2a\leq b$, and so by (\ref{ineqI1c}) we obtain that $2a=b$. Thus
(up to the action of $GL(3,\Z)$) we have that $a=1, b=2, c=3$. A
quick check confirms that this is indeed Fano.
\end{proof}
\begin{lem}\label{lemCaseI2}
If $P$ is a minimal Fano polytope with vertices
$\set{x_1,x_2,x_3,x_4,x_5}$ such that $\set{x_2,x_3,x_4,x_5}$ are
coplanar and give the vertices of a parallelogram, then $P$ is
equal (up to the action of $GL(3,\Z)$) to
$$\begin{pmatrix}1&0&-1&1&0\\0&1&-1&1&0\\0&0&0&1&-1\end{pmatrix}.$$
\end{lem}
\begin{proof}
Since $P$ does not contain a Fano tetrahedron it must be that
opposite corners of the parallelogram, along with $x_1$, give us a
Fano triangle. Thus we can (by virtue of Lemma~\ref{emptyTriEqu})
write $P$ in the form
$$\begin{pmatrix}1&0&-1&a+1&-a\\0&1&-1&b+1&-b\\0&0&0&c&-c\\\end{pmatrix}$$
where $0<a+1\leq b+1\leq c$.

Consider the point $-e_3$. Either $a=0,b=0,c=1$, which gives the
Fano polytope in the statement, or $-e_3$ lies outside $P$. Thus
$-e_3$ lies on the opposite side to the origin of the plane
defined by $\set{e_1,e_2,\coord{-a,-b,-c}}$. This plane intersects
the $z$-axis at the point $\coord{0,0,\frac{c}{a+b+1}}$. Thus we
have that $-c>-a-b-1$ and so
\begin{equation}\label{ineq1l2}
c\leq a+b
\end{equation}

Now let $x'=e_1+e_2+e_3$. Either $a=0,b=0,c=1$, which gives the
Fano polytope in the statement, or $x'$ lies outside $P$. Thus
$x'$ lies on the opposite side to the origin of the plane defined
by $\set{e_1,e_2,\coord{a+1,b+1,c}}$. Thus plane intersects the
line through the origin and $x'$ at the point $\coord{k,k,k}$,
where $k:=\frac{c}{2c-a-b-1}$. Thus we see that $c<2c-a-b-1$ and
so
\begin{equation}\label{ineq2l2}
c>a+b+1
\end{equation}

Now suppose both $-e_3$ and $x'$ lie outside $P$. But then both
inequalities (\ref{ineq1l2}) and (\ref{ineq2l2}) must be
satisfied, which is impossible.
\end{proof}
Combining the results of Table \ref{tetTable} and Lemmas
\ref{lemCaseII}-\ref{lemCaseI2} we obtain Table \ref{minTab}.
\begin{table}
\centering
\begin{tabular}{cc}
\begin{tabular}{|c|c|}
\hline Comments&Vertices\\
\hline
\begin{tabular}{c}4 Vertices\\Simplicial\end{tabular}&$\left(\begin{array}{cccc}
1&0&1&-2\\0&1&-3&2\\0&0&5&-5\end{array}\right)$\\
\hline\begin{tabular}{c}4
Vertices\\Simplicial\end{tabular}&$\left(\begin{array}{cccc}
1&0&1&-1\\0&1&-2&1\\0&0&7&-5\end{array}\right)$\\
\hline\begin{tabular}{c}4
Vertices\\Simplicial\end{tabular}&$\left(\begin{array}{cccc}
1&0&-2&1\\0&1&2&-2\\0&0&7&-5\end{array}\right)$\\
\hline\begin{tabular}{c}4
Vertices\\Simplicial\end{tabular}&$\left(\begin{array}{cccc}
1&0&0&-1\\0&1&0&-1\\0&0&1&-1\end{array}\right)$\\
\hline\begin{tabular}{c}4
Vertices\\Simplicial\end{tabular}&$\left(\begin{array}{cccc}
1&0&1&-1\\0&1&-2&1\\0&0&5&-4\end{array}\right)$\\
\hline\begin{tabular}{c}4
Vertices\\Simplicial\end{tabular}&$\left(\begin{array}{cccc}
1&0&-2&1\\0&1&1&-1\\0&0&5&-3\end{array}\right)$\\
\hline\begin{tabular}{c}4
Vertices\\Simplicial\end{tabular}&$\left(\begin{array}{cccc}
1&0&1&-1\\0&1&1&-1\\0&0&2&-1\end{array}\right)$\\
\hline
\end{tabular}
&\begin{tabular}{|c|c|}
\hline Comments&Vertices\\
\hline
\begin{tabular}{c}4 Vertices\\Simplicial\end{tabular}&$\left(\begin{array}{cccc}
1&0&1&-1\\0&1&-2&1\\0&0&3&-2\end{array}\right)$\\
\hline\begin{tabular}{c}5
Vertices\\Simplicial\end{tabular}&$\left(\begin{array}{ccccc}
1&0&1&-1&-1\\0&1&2&-1&-2\\0&0&3&0&-3\end{array}\right)$\\
\hline\begin{tabular}{c}5
Vertices\\Simplicial\end{tabular}&$\left(\begin{array}{ccccc}
1&0&0&-1&0\\0&1&0&-1&0\\0&0&1&0&-1\end{array}\right)$\\
\hline\begin{tabular}{c}5
Vertices\end{tabular}&$\left(\begin{array}{ccccc}
1&0&1&-1&0\\0&1&1&-1&0\\0&0&1&0&-1\end{array}\right)$\\
\hline\begin{tabular}{c}6
Vertices\\Simplicial\end{tabular}&$\left(\begin{array}{cccccc}
1&0&1&-1&0&-1\\0&1&1&0&-1&-1\\0&0&2&0&0&-2\end{array}\right)$\\
\hline\begin{tabular}{c}6
Vertices\\Simplicial\end{tabular}&$\left(\begin{array}{cccccc}
1&0&0&-1&0&0\\0&1&0&0&-1&0\\0&0&1&0&0&-1\end{array}\right)$\\
\hline
\end{tabular}
\end{tabular}
\caption{The vertices of the minimal Fano polytopes, up to the
action of $GL(3,\Z)$}\label{minTab}
\end{table}
\section{Classifying all Fano Polytopes}\label{finalClassAll}
Given any Fano polytope $P$ with vertices $\set{x_1,\ldots,x_k}$
we make the following definition (c.f. Definition
\ref{minimalDef}):
\begin{defn}\label{maximalDef}
We say $P$ is \emph{maximal} if, for all
$x_{k+1}\in\Z^3\setminus\set{x_1,\ldots,x_k}$, the polytope $P''$
given by the vertices $\set{x_1,\ldots,x_k,x_{k+1}}$ is not Fano.
\end{defn}
\begin{table}[ht]
\centering
\begin{tabular}{|c|l|}
\hline Comments&Vertices\\
\hline
\begin{tabular}{c}8 Vertices\\Simplicial\end{tabular}&$\left(\begin{array}{cccccccc}
1&0&0&-1&-1&0&-1&3\\0&1&0&-1&0&-1&1&-2\\0&0&1&-1&0&-1&2&-1\end{array}\right)$\\
\hline\begin{tabular}{c}8
Vertices\\Simplicial\end{tabular}&$\left(\begin{array}{cccccccc}
1&0&0&-1&-1&1&-2&3\\0&1&0&-1&0&-1&-1&-2\\0&0&1&-1&0&0&-1&-1\end{array}\right)$\\
\hline\begin{tabular}{c}8
Vertices\\Simplicial\end{tabular}&$\left(\begin{array}{cccccccc}
1&0&1&-2&-1&1&0&0\\0&1&-3&2&1&-1&-1&1\\0&0&5&-5&-2&2&1&-1\end{array}\right)$\\
\hline\begin{tabular}{c}9
Vertices\end{tabular}&$\left(\begin{array}{ccccccccc}
1&0&0&-1&-1&0&1&-1&-2\\0&1&0&-1&0&-1&1&-2&1\\0&0&1&-1&0&0&0&-1&-1\end{array}\right)$\\
\hline\begin{tabular}{c}10
Vertices\end{tabular}&$\left(\begin{array}{cccccccccc}
1&0&0&-1&-1&0&0&-1&0&-1\\0&1&0&-1&0&-1&0&1&-1&2\\0&0&1&-1&0&0&-1&0&1&1\end{array}\right)$\\
\hline\begin{tabular}{c}10
Vertices\end{tabular}&$\left(\begin{array}{cccccccccc}
1&0&0&-1&-1&0&0&-1&0&-1\\0&1&0&-1&0&-1&0&1&1&-2\\0&0&1&-1&0&0&-1&0&-1&-1\end{array}\right)$\\
\hline\begin{tabular}{c}11
Vertices\end{tabular}&$\left(\begin{array}{ccccccccccc}
1&0&0&-1&-1&0&0&1&-1&0&1\\0&1&0&-1&0&-1&0&1&1&-1&0\\0&0&1&-1&0&0&-1&1&0&1&-1\end{array}\right)$\\
\hline\begin{tabular}{c}12
Vertices\end{tabular}&$\left(\begin{array}{cccccccccccc}
1&0&0&-1&-1&0&0&1&-1&1&-1&0\\0&1&0&-1&0&-1&0&1&-1&1&0&1\\0&0&1&-1&0&0&-1&1&0&0&1&-1\end{array}\right)$\\
\hline\begin{tabular}{c}14
Vertices\end{tabular}&$\left(\begin{array}{cccccccccccccc}
1&0&0&-1&-1&0&0&1&-1&1&-1&1&0&0\\0&1&0&-1&0&-1&0&1&-1&1&0&0&1&-1\\0&0&1&-1&0&0&-1&1&0&0&-1&1&1&-1\end{array}\right)$\\
\hline
\end{tabular}
\caption{The vertices of the maximal Fano polytopes, up to the
action of $GL(3,\Z)$}\label{maxTab}
\end{table}
We will also make the following non-standard definition:
\begin{defn}
Let $P=\sconv{x_1,\ldots,x_k}$ and $P''$ be Fano polytopes and let
$x_{k+1}\in\Z^3$ be a point such that, up to the action of
$GL(3,\Z)$, $P''=\sconv{x_1,\ldots,x_k,x_{k+1}}$. Then we say that
$P$ is the \emph{parent} of $P''$, and that $P''$ is the
\emph{child} of $P$.
\end{defn}
Clearly a polytope $P$ is minimal if and only if it has no
parents, and is maximal if and only if it has no children.

Let $P$ be any Fano polytope. Then the following results are
immediate:
\begin{enumerate}
 \item[(i)]Any Fano polytope can be obtained from a (not
 necessarily unique) minimal Fano polytope by consecutive
 addition of vertices.
 \item[(ii)]The number of possible vertices that can be added to
 $P$ to create a Fano polytope $P''$ is finite. For suppose
 $P$ has vertices $\set{x_1,\ldots,x_n}$ and the vertex $x_{n+1}$
 is to be added. Then the line through $x_{n+1}$ and the origin,
 extended in the direction away from $x_{n+1}$, crosses $\partial P$
 at either a vertex point, an edge, or a face. The first possibility gives us that
 $x_{n+1}=-x_i$ for some $i\in\set{1,\ldots,n}$. The second
 possibility tells us that $\set{x_i,x_j,x_{n+1}}$ is an
 Fano triangle for some distinct $i,j\in\set{1,\ldots,n}$,
 and hence that $x_{n+1}=-x_i-x_j$. The final possibility
 corresponds to $\set{x_i,x_j,x_k,x_{n+1}}$ being a Fano
 tetrahedron for some distinct $i,j,k\in\set{1,\ldots,n}$, and so
 $\lambda_{\sigma1}x_i+\lambda_{\sigma2}x_j+\lambda_{\sigma3}x_k+\lambda_{\sigma4}x_{n+1}=0$ for
 some $\coord{\lambda_1,\ldots,\lambda_4}$ in Proposition
 \ref{mainProp} and some $\sigma\in S_4$.
 \item[(iii)]If $\set{x_1,\ldots,x_n}$ are the vertices of $P$,
 and the Fano polytope $P''$ is created by adding the vertex
 $x_{n+1}$, then
 $$P''\setminus P\subset\bigcup_{i,j}\sconv{0,x_i,x_j,x_{n+1}}.$$
\end{enumerate}

Using these results and our list of minimal Fano polytopes, it is
a relatively straightforward task to write a recursive function to
allow a computer to calculate all the Fano polytopes up to the
action of $GL(3,\Z)$. In particular, (ii) asserts that the
calculation will terminate, since the list is finite; a stronger
finiteness result to include $\varepsilon$-logcanonical toric Fano
varieties ($0<\varepsilon\leq 1$) can be found
in~\cite{BB92,Bor00}.

The source code for such a function is available on the internet
at
\begin{center}
\texttt{http://www.maths.bath.ac.uk/$\sim$mapamk/code/Polytope\underline{
}Classify.c}.
\end{center}

Using this code a complete classification was obtained in under
$20$ minutes on an average personal computer. This list, along
with a table giving the parents and children of each Fano
polytope, is available on the internet at
\begin{center}
\texttt{http://www.maths.bath.ac.uk/$\sim$mapamk/pdf/Fano\underline{
}List.pdf}
 (or \texttt{.ps}).
\end{center}

The maximal polytopes are reproduced in Table \ref{maxTab}, and a
summaries of the results are given in Theorem~\ref{mainSummary}
and in Table \ref{3DNumTable}.
\begin{table}[ht]
\centering
\begin{tabular}{|r|ccccccccccc|}
\hline Vertices&4&5&6&7&8&9&10&11&12&13&14\\
\hline\hline Polytopes&8&38&95&144&151&107&59&21&8&2&1\\
\hline Simplicial&8&35&75&74&35&5&1&0&0&0&0\\
\hline Minimal&8&3&2&0&0&0&0&0&0&0&0\\
\hline Maximal&0&0&0&0&3&1&2&1&1&0&1\\
\hline
\end{tabular}
\caption{The number of Fano polytopes in
$\mathbb{Z}^3$.}\label{3DNumTable}
\end{table}
\bibliographystyle{amsalpha}
\providecommand{\bysame}{\leavevmode\hbox
to3em{\hrulefill}\thinspace}
\providecommand{\MR}{\relax\ifhmode\unskip\space\fi MR }
\providecommand{\MRhref}[2]{%
  \href{http://www.ams.org/mathscinet-getitem?mr=#1}{#2}
} \providecommand{\href}[2]{#2}

\end{document}